\documentclass[12pt,reqno]{amsart}
\usepackage[utf8]{inputenc}
\usepackage[T1]{fontenc}
\usepackage{amsxtra,amssymb,amsthm,amsmath,amscd,mathrsfs, epsfig, eufrak}
\usepackage{amscd, amsmath, mathrsfs, amssymb, amsthm, amsxtra, bbding, epsfig, eucal, eufrak, graphicx, latexsym, mathrsfs, url, color, mathbbol, bbold}
\usepackage{mathtools}
\usepackage[all]{xy}
\usepackage{bbding}
\usepackage{float}
\usepackage{graphicx}
\usepackage{latexsym}
\usepackage{epsfig,epstopdf,color}
\usepackage{bbold}
\usepackage{bbding}
\usepackage{oldgerm}
\usepackage[ngerman,french,english]{babel}
\usepackage{titletoc}

\usepackage[colorlinks,linkcolor=blue,anchorcolor=blue,citecolor=blue]{hyperref}

\theoremstyle{plain}
\newtheorem{Theorem}{Theorem}
\newtheorem*{Theorem*}{Theorem}
\newtheorem{Lemma}{Lemma}[section]
\newtheorem*{hypothesis*}{Hypothesis}
\newtheorem{corollary}{Corollary}

\newtheorem*{corollaire*}{Corollaire 2*}

\newtheorem*{conjecture*}{Conjecture}

\theoremstyle{definition}

\newtheorem{remark}{Remark}
\numberwithin{equation}{section}

\usepackage{tikz}
\usepackage[normalem]{ulem}
\usepackage{soul}
\usepackage{color}
\setstcolor{red}

\def\leq{\leqslant}
\def\geq{\geqslant}

\def\d{{\rm{d}}}
\def\e{{\rm{e}}}

\allowdisplaybreaks

\begin{document}

\title[Averages of coefficients of a class of degree 3 $L$-functions]
{Averages of coefficients of a class of degree 3 $L$-functions}
\date{\today}
\author{Bingrong Huang, Yongxiao Lin, and Zhiwei Wang}

\address{%
Bingrong Huang \newline
\indent$\prescript{}{}{}$Data Science Institute and School of Mathematics\\
Shandong University\\
Jinan\\
250100 Shandong\\
P. R. China
}
\email{brhuang@sdu.edu.cn}
\address{%
Yongxiao Lin\newline
\indent$\prescript{}{}{}$EPFL SB MATHGEOM TAN\\
Station 8, CH-1015, Lausanne\\
Switzerland}
\email{yongxiao.lin@epfl.ch}
\address{%
Zhiwei Wang\newline
\indent$\prescript{}{}{}$School of Mathematics\\
Shandong University\\
Jinan\\
250100 Shandong\\
P. R. China}
\email{zhiwei.wang@sdu.edu.cn}

\selectlanguage{english}
\begin{abstract}
  In this note, we give a detailed proof of an asymptotic for averages of coefficients of a class of degree three $L$-functions which can be factorized as a product of a degree one and a degree two $L$-functions. We emphasize that we can break the $1/2$-barrier in the error term, and we get an explicit exponent.
\end{abstract}

\thanks{B. H. was supported by the Young Taishan Scholars Program of Shandong Province (Grant No. tsqn201909046) and Qilu Young Scholar Program of Shandong University.
Y. L. was partially supported by a DFG-SNF lead agency program grant (Grant  200020L\_175755).
Z. W. was supported by China Postdoctoral Science Foundation (No. 2019M652354) and  National Natural Science Foundation of China (No. 11901348)}

\vglue -1,5mm
\maketitle
{\footnotesize

\dottedcontents{section}[1.16cm]{}{1.8em}{5pt}
\dottedcontents{subsection}[2.00cm]{}{2.7em}{5pt}
}

\section{Introduction}

In number theory, one of the most important and central problem is to give an asymptotic formula as accurate as possible for the sum
\begin{align*}
A(X)=\sum_{n\leq X}a(n)
\end{align*}
where $a(n)$ is an arithmetic function. In many cases, this problem can be transformed to the study of the following generating series
\begin{align*}
  L(s)=\sum_{n=1}^{\infty}\frac{a(n)}{n^s}
\end{align*}
which satisfies valuable analytic properties. Here we just list two classical examples.
\begin{itemize}
\item
$a(n)=\Lambda(n)$ is the von Mangoldt function, then this is the well-known prime number theorem.

\vskip 2mm

\item
$a(n)=\tau_k(n), k\geq 2$ is the divisor function, then this is the divisor problem. One may see for example the work of Atkinson \cite{Atkinson41} for the case $k=3$.
\end{itemize}

\vskip 0.6mm

In this paper, we are interested in investigating the case where $a(n)=\lambda_F(n)$, the coefficients of a degree $d$ automorphic $L$-function (see Section \ref{L-functions}). We set
\begin{align*}
A(X, F)=\sum_{n\leq X}\lambda_F(n),\qquad L(s, F)=\sum_{n\geq 1}\frac{\lambda_F(n)}{n^{s}}.
\end{align*}

\vskip 1,3mm

For $A(X, F)$, one may expect that
\begin{align}\label{asymptotic}
A(X, F)=\mathop{\text{Res}}\limits_{s=1}\frac{L(s, F)}{s}X+O(X^{\theta})
\end{align}
with $\theta$ as small as possible.
In particular, this would imply
\begin{align}\label{short-interval}
\sum_{X<n\leq X+Y}\lambda_F(n)=\mathop{\text{Res}}\limits_{s=1}\frac{L(s, F)}{s}Y+ o(Y)
\end{align}
for $Y>X^{\theta+\varepsilon}$, for any fixed $\varepsilon>0$.

Under the generalized Riemann hypothesis (GRH), we may take $\theta=1/2$. To surpass the GRH with some $\theta<1/2$ is an interesting and challenging problem. For general series of degree 2, this is already achieved with $\theta=1/3$. Friedlander and Iwaniec \cite{FI2005} proved a quite general result that for $d\geq 2$, \eqref{asymptotic} holds for $\theta=\frac{d-1}{d+1}$ ($\geq \frac{1}{2}$ for $d\geq 3$) under some local assumptions.

\vskip 1,7mm

For the case of degree ``1+1+1'' when $L(s, F)$ factors completely into degree 1 $L$-functions, it was showed that we can take $\theta=\frac{37}{75}<\frac{1}{2}$ in \eqref{asymptotic} (see e.g. Friedlander and Iwaniec \cite{FI2005}). They also announced that the other case of degree ``1+2'' can also be successfully treated with additional arguments. For $d=6$, recently together with Sun, the second named author \cite{LinSun} was able to beat the bound $\frac{6-1}{6+1}$ (under Ramanujan--Selberg) for a class of $\rm GL_3\times \rm GL_2$ automorphic forms, namely, for $F=f\otimes\phi$, Rankin--Selberg convolution of a $\rm GL_3$ cusp form $f$ and a $\rm GL_2$ cusp form $\phi$.  Subsequently, the first named author \cite{Huang2020} was able to make further progress on another important case, by improving the error term $O(X^{3/5})$ for $F=\phi\times \phi$ for the case $d=4$. Specially, he remarked that the method in \cite{Huang2020} is applicable to the degree ``1+2'' case. Here, crucial ingredients in these works are factorazation of the coefficients: in Friedlander--Iwaniec's case $\lambda_F(n)=\sum_{n_1n_2n_3=n}\chi_1(n_1)\chi_2(n_2)\chi_3(n_3)$; in Lin--Sun's case $\lambda_F(n)=\sum_{m_1^2m_2=n}\lambda_f(m_1,m_2)\lambda_\phi (m_2)$, while in Huang's case $\lambda_F(n)=\sum_{\ell m=n}\lambda_{\text{sym}^2}(m)$.

In this paper, we shall consider $A(X, F)$ for another case of degree 3, when $F=1\boxplus \phi$; that is
\begin{align*}
L(s, F)=L(s, 1\boxplus \phi)=\zeta(s)L(s, \phi).
\end{align*}
Here $\phi$ is a holomorphic cusp form. We give a first detailed proof of \eqref{asymptotic} with $\theta<\frac{1}{2}$, by combining the methods of Friedlander--Iwaniec in \cite{FI2005} and the first named author in \cite{Huang2020}.

\begin{Theorem}\label{thm1}
With the notation as above. The asymptotic formula
\begin{align*}
A(X, 1\boxplus \phi)=L(1, \phi)X+O(X^{1/2-\delta+\varepsilon})
\end{align*}
holds for $\delta= 4/739 \approx 0.00541...$ and any fixed $\varepsilon>0$.
\end{Theorem}

An immediate consequence of this theorem is the following asymptotic formula for averages of those coefficients in rather short intervals.
\begin{corollary}
With the notation as above. We have
\begin{align*}
\sum_{X<n\leq X+Y}\lambda_{1\boxplus \phi}(n)=L(1, \phi)Y+O(Y^{1-\varepsilon})
\end{align*}
as long as $Y\gg X^{1/2-\delta}$, for any $\delta< 4/739 \approx 0.00541...$.
\end{corollary}

The functions $n\mapsto \tau_3(n)$ and $n\mapsto \lambda_{1\boxplus\phi}(n)$ are the Hecke eigenvalues of certain non-cuspidal automorphic representation of $\operatorname{GL}_{3,\mathbb{Q}}$, that is, the isobaric representation $1\boxplus1\boxplus 1$ and $1\boxplus \pi_\phi$.
The methods of \cite{Atkinson41,FI2005,Huang2020} and of the present paper can be generalized straightforwardly to prove that the average of the $n$-th Hecke eigenvalue function of any fixed non-cuspidal automorphic representation of $\operatorname{GL}_{3,\mathbb{Q}}$ has exponent $<1/2$ in the error term.
Extending this further to cuspidal $\operatorname{GL}_{3,\mathbb{Q}}$-representations is an interesting and challenging question.

\begin{remark}
  Under the generalized Ramanujan conjecture, it is clear that the same result can be proved for $\phi$ being a Hecke--Maass cusp form by our method.
\end{remark}

\begin{remark}
(1) There is an arithmetic analogue of the questions \eqref{asymptotic} and \eqref{short-interval}. Indeed, one can study the distribution of $\lambda_F(n)$ by also asking how uniform the coefficients $\lambda_F(n)$ are distributed, when we vary $n$ among arithmetic progressions $am+q$. Namely, do we have
\begin{equation}\label{level-of-distribution}
\mathop{\sum_{n\leq X}}_{n\equiv a \bmod q}\lambda_F(n)-\frac{1}{\varphi(q)}\mathop{\sum_{n\leq X}}_{(n,q)=1}\lambda_F(n)\ll_A \frac{X}{q}\left(\log X\right)^{-A}
\end{equation}
for $q\leq X^{\vartheta}$, with $\vartheta$ as large as possible? Here $\varphi$ is the Euler function and the exponent $\vartheta$ is called the level of distribution. It is predicted that one can take $\vartheta=1-\varepsilon$. In general, for $F$ an $\rm GL_d$ automorphic form, this can be studied by detecting the condition $n\equiv a \bmod q$ using additive characters and then applying the Vorono\"{i} summation, which would transform the sum in question to a dual sum of length roughly $q^d/X$ and to the dual Hecke eigenvalue $\overline{\lambda_F(n)}$, but twisted by a $(d-1)$-dimensional hyper-Kloosterman sum $\text{Kl}_d(an;q)$, upon which an application of Deligne's estimate will produce a level of distribution $\frac{2}{d+1}$ for the sum above. Such an exponent is regarded as the trivial level of distribution for $\lambda_F(n)$. To beat this exponent for various $\lambda_F$ is an active area of research.
For example, results breaking the $1/2$-barrier for the level $\vartheta$ which is also the same barrier of the well-known Bombieri--Vinogradov theorem, are known for several special cases:
\begin{itemize}
\item
when $\lambda_F(n)=\tau(n)$, we can take $\vartheta=2/3-\varepsilon$, a classical result of Selberg which remains the best to date;

\vskip 2mm

\item
 when $\lambda_F(n)=\tau_3(n)$, one can take $\vartheta=1/2+1/230-\varepsilon$, corresponding to the groundbreaking result of Friedlander--Iwaniec \cite{FI1985}.
 \vskip 2mm

\item
when $\lambda_F(n)=\lambda_{1\boxplus\phi}(n)$, one has $\vartheta=1/2+1/102-\varepsilon$, which was proved by Kowalski--Michel--Sawin \cite{KMS}, building on their breakthrough on estimates for bilinear forms in hyper-Kloosterman sums.

\end{itemize}

\par
(2) Theorem \ref{thm1} (and its corollary) is exactly an Archimedian analogue of \eqref{level-of-distribution} when $F=1\boxplus\phi$, proved by Kowalski--Michel--Sawin.
\end{remark}

Throughout the paper, $\varepsilon$ is an arbitrarily small positive number;
all of them may be different at each occurrence.

\vskip 7mm
\section{Preliminaries}
\subsection{$L$-functions}\label{L-functions}
In general, the $L$-function $L(s, F)$ satisfies the following conditions:

(1).\, We have the Euler product of degree $d$
\begin{align*}
L(s, F)=\sum_{n\geq 1}\frac{\lambda_F(n)}{n^{s}}=\prod_{p}\prod_{j=1}^{d}\Big(1-\frac{\alpha_j(p)}{p^s}\Big)^{-1}
\end{align*}
with $\lambda_f(1)=1$. 
The series and Euler products are
absolutely convergent for Re$\,s > 1$. The sequence $\{\lambda_F(n)\}_{n\geq 1}$ are called coefficients of $L(s, F)$, and the $\alpha_j(p), 1\leq j\leq d$ satisfying $|\alpha_j(p)|<p$ for all $p$, are called the local parameters of $L(s, F)$ at $p$.

(2).\, We have the gamma factor defined by
\begin{align*}
\gamma(s, F)=\pi^{-ds/2}\prod_{j=1}^{d}\Gamma\Big(\frac{s+\kappa_j}{2}\Big)
\end{align*}
where the numbers $\kappa_j\in \mathbb{C}$ satisfying Re$\,\kappa_j>-1$, are called the local parameters of $L(s, F)$ at infinity.

(3).\, There exists an integer $q(F)\geq 1$, called the conductor of $L(s, F)$, satisfying $\alpha_j(p)\neq 0$ for $1\leq j\leq d$, $p\nmid q(F)$ such that we have the functional equation
\begin{align*}
\Lambda(s, F)=\varepsilon(F)\Lambda(1-s, \bar{F})
\end{align*}
where $\bar{F}$ is the dual form of $f$ for which $\lambda_{\bar{F}}(n)=\overline{\lambda_F(n)}$,
$\gamma(s, \bar{F})=\gamma(s, F)$, $q(\bar{F})=q(F)$, and $\varepsilon(F)$ is the root number of $L(s, F)$ satisfying $|\varepsilon(F)|=1$. Here $\Lambda(s, F)$ is called the complete $L$-function defined by
\begin{align*}
\Lambda(s, F)=q(F)^{s/2}\gamma(s, F)L(s, F).
\end{align*}
For more details, we refer the reader to the book of Iwaniec--Kowalski \cite[Chapter 5.1]{IwaniecKowalski2004analytic}.

\subsection{Functional equation}
\begin{Lemma}\label{lemma}
For {\rm Re}$\,s>1$, we have
\begin{align*}
L(1-s,\, 1\boxplus \phi)=w\gamma(s)L(s,\, 1\boxplus\phi)
\end{align*}
where $w=i^k$ and $\gamma(s)$ satisfies
\begin{align*}
\gamma(\sigma-it)=\overline{\omega}\, (Qt)^{3(\sigma-1/2)}
\Big(\frac{\e}{Qt}\Big)^{3it}\Big\{1+O\Big(\frac{1}{t}\Big)\Big\}
\end{align*}
for $\sigma > 1/2$, $t > 1$, and
$$
Q=\frac{1}{2\pi}, \qquad \omega=\e\Big(\frac{2k-3}{8}\Big),\ \, where \ \, \e(u)=\e^{2\pi i u}.
$$
\end{Lemma}

\vskip 3mm

\begin{proof}


First by the functional equation
$$
\Lambda(s,\, 1\boxplus \phi)=\pi^{-\frac{s}{2}}\Gamma\big(\frac{s}{2}\big)\zeta(s)\cdot
(2\pi)^{-s}\Gamma\Big(s+\frac{k-1}{2}\Big)L(s,\, \phi)=i^k\Lambda(1-s,\, 1\boxplus \phi),
$$
we can write the functional equation of the Dirichlet generating series for the coefficient $\lambda_{1\boxplus \phi}(n)$ as follows:
\begin{equation}\label{L(1-s)}
    \begin{aligned}
L(1-s,\, 1\boxplus \phi)&=i^k\frac{\pi^{-\frac{s}{2}}\Gamma\big(\frac{s}{2}\big)(2\pi)^{-s}\Gamma\Big(s+\frac{k-1}{2}\Big)}
{\pi^{-\frac{1-s}{2}}\Gamma\big(\frac{1-s}{2}\big)(2\pi)^{-1+s}\Gamma\Big(1-s+\frac{k-1}{2}\Big)} L(s,\, 1\boxplus\phi)
\\\noalign{\vskip 2,8mm}
&=w\gamma(s)L(s,\, 1\boxplus\phi)
\end{aligned}
\end{equation}
where
\begin{align}\label{w+gamma}
w=i^k,\qquad \gamma(s)=\pi^{\frac{1}{2}-s}(2\pi)^{1-2s}\, \frac{\Gamma\big(\frac{s}{2}\big)}{\Gamma\big(\frac{1-s}{2}\big)}\cdot
\frac{\Gamma\Big(s+\frac{k-1}{2}\Big)}{\Gamma\Big(1-s+\frac{k-1}{2}\Big)}
\end{align}
Next, we shall write $\gamma(s)$ in the form
\begin{align}\label{gamma(s)}
\gamma(s)=(\pi^{-m}D)^{s-\frac{1}{2}}\prod_{j=1}^{m}
\Gamma\Big(\frac{s+\kappa_j}{2}\Big)\Gamma\Big(\frac{1-s+\kappa_j}{2}\Big)^{-1}
\end{align}
where $D$ is the conductor, and the “spectral parameters” $\kappa_j$
are complex numbers having Re$\,\kappa_j > 0$. And $k=\sum_{j=1}^{m}\kappa_j$ is  the weight.

Note that $\Gamma(z)\Gamma(z+\frac{1}{2})=2^{1-2z}\pi^{\frac{1}{2}}\Gamma(2z)$. Taking $z=\frac{s+\frac{k-1}{2}}{2}$, we obtain
$$
\Gamma\Big(s+\frac{k-1}{2}\Big)=2^{s+\frac{k-3}{2}}\pi^{-\frac{1}{2}}\Gamma\Big(\frac{s+\frac{k-1}{2}}{2}\Big)
\Gamma\Big(\frac{s+\frac{k+1}{2}}{2}\Big).
$$
With the notation of $\gamma(s)$ in \eqref{w+gamma}, we derive
\begin{equation}\label{gamma(s)-2}
    \begin{aligned}
\gamma(s)&= \pi^{\frac{1}{2}-s}(2\pi)^{1-2s}\, \frac{\Gamma\big(\frac{s}{2}\big)}{\Gamma\big(\frac{1-s}{2}\big)}\cdot
\frac{2^{s+\frac{k-3}{2}}\pi^{-\frac{1}{2}}\Gamma\Big(\frac{s+\frac{k-1}{2}}{2}\Big)
\Gamma\Big(\frac{s+\frac{k+1}{2}}{2}\Big)}{2^{1-s+\frac{k-3}{2}}\pi^{-\frac{1}{2}}\Gamma\Big(\frac{1-s+\frac{k-1}{2}}{2}\Big)
\Gamma\Big(\frac{1-s+\frac{k+1}{2}}{2}\Big)}
\\\noalign{\vskip 2,8mm}
&=(\pi^{-3})^{s-\frac{1}{2}}\prod_{j=1}^{3}\Gamma\Big(\frac{s+\kappa_j}{2}\Big)\Gamma\Big(\frac{1-s+\kappa_j}{2}\Big)^{-1},
\end{aligned}
\end{equation}
which is of the form \eqref{gamma(s)} with $m=3$, $D=1$ and $\kappa_1=0$, $\kappa_2=\frac{k-1}{2}$, $\kappa_2=\frac{k+1}{2}$.
Hence, by the argument of Friedlander--Iwaniec \cite[Section 1]{FI2005}, we complete the proof of Lemma \ref{lemma} with $Q=(2\pi)^{-1}$, $\omega=\e(\frac{2k-3}{8})$.
\end{proof}

On taking $s=1+\varepsilon-it$, Lemma \ref{lemma} yields
\begin{equation}\label{functional eq}
L(-\varepsilon+it,\, 1\boxplus \phi)
=w\,\overline{\omega}\, \Big(\frac{t}{2\pi}\Big)^{3(\frac{1}{2}+\varepsilon)}\Big(\frac{t}{2\pi \e}\Big)^{-3it}
L(1+\varepsilon-it,\, 1\boxplus\phi)\Big\{1+O\Big(\frac{1}{t}\Big)\Big\}.
\end{equation}

\vskip 7mm
\section{Proof of Theorem \ref{thm1}}
In this section, we consider the sum $A(X,\, 1\boxplus \phi)=\sum_{n\leq X}\lambda_{1\boxplus \phi}(n)$, that is we have
$$L(s,\, 1\boxplus \phi)=\zeta(s)L(s,\, \phi)=\zeta(s)\sum_{n=1}^{\infty}\frac{\lambda_{\phi}(n)}{n^s}$$
for Re$(s)>1$, where $\phi$ is a $\text{GL}_2$ automorphic form and $\lambda_{\phi}(n)$ is defined by
$$
\phi(z)=\sum_{n\geq 1}\lambda_{\phi}(n)n^{\frac{k-1}{2}}\text{e}(nz),\quad \text{e}(z)=\text{e}^{2\pi iz}.
$$
Assume $\phi$ is holomorphic, i.e., $\phi\in H_k(\text{SL}_2(\mathbb{Z}))$, an orthogonal basis of the space of holomorphic cusp forms of weight $k$ and level $1$, with $k\geq 12$ even.
We first approximate $A(X,\, 1\boxplus \phi)$ by a smooth sum. Let
$$
Y=X^{1/2-\delta},\quad \text{for some}\ \delta\in(0,\, 1/10).
$$
Let $W$ be a smooth function with support supp$\,W\in[1/2-Y/X,\, 1+Y/X]$ such that $W(u)=1, u\in[1/2,\, 1]$ and $W(u)\in[0,\, 1], u\in [1/2-Y/X,\, 1/2] \cup [1,\, 1+Y/X]$, and $W^{(k)}(u)\ll (X/Y)^k$. Therefore we have the approximating forumla
\begin{equation}\label{Approx}
    \begin{aligned}
\sum_{X/2<n\leq X}\lambda_{1\boxplus \phi}(n)=&\sum_{X/2-Y<n<X+Y}\lambda_{1\boxplus \phi}(n)W\Big(\frac{n}{X}\Big)
\\\noalign{\vskip 2,3mm}
\quad &+O\bigg(\sum_{X/2-Y<n<X/2}|\lambda_{1\boxplus \phi}(n)|+\sum_{X<n<X+Y}|\lambda_{1\boxplus \phi}(n)|\bigg)
\\\noalign{\vskip 2,8mm}
=&\sum_{n\geq 1}\lambda_{1\boxplus \phi}(n)W\Big(\frac{n}{X}\Big)+O\big(X^{1/2-\delta+\varepsilon}\big)
\end{aligned}
\end{equation}
where we have used Deligne's bound $\lambda_{1\boxplus \phi}(n)\ll \sum_{\ell m=n}  m^{\varepsilon}  \ll n^{\varepsilon}$ when $\phi$ is holomorphic for the error terms.

Next we only need to show
\begin{align}\label{Aim}
\sum_{n\geq 1}\lambda_{1\boxplus \phi}(n)W\Big(\frac{n}{X}\Big)=L(1,\, \phi)\tilde{W}(1)X+O\big(X^{1/2-\delta+\varepsilon}\big)
\end{align}
where $\tilde{W}(s)=\int_{0}^{\infty}W(x)x^{s-1}\d x$ is the Mellin transform of $W$. Since by inserting \eqref{Aim} into \eqref{Approx} and then by iteration, we get
\begin{align*}
\sum_{n\leq X}\lambda_{1\boxplus \phi}(n)=2L(1,\, \phi)\tilde{W}(1)X+O\big(X^{1/2-\delta+\varepsilon}\big).
\end{align*}
Then Theorem \ref{thm1} follows immediately from the estimate $\tilde{W}(1)=1/2+O(Y/X)$.

Now we estimate the sum $\sum_{n\geq 1}\lambda_{1\boxplus \phi}(n)W(\frac{n}{X})$ in \eqref{Aim}. By the inverse Mellin transform
$$
W(u)=\frac{1}{2\pi i}\int_{(2)}\tilde{W}(s)u^{-s}\d s,
$$
we get
\begin{align*}
\sum_{n\geq 1}\lambda_{1\boxplus \phi}(n)W\Big(\frac{n}{X}\Big)=\frac{1}{2\pi i}\int_{(2)}\tilde{W}(s)L(s,\, 1\boxplus \phi)X^{s}\d s.
\end{align*}
We then move the integration to the parallel segment with Re$\,s=\sigma=-\varepsilon$. We pass pole at $s=1$ with residue
Res$_{s=1}L(s,\, 1\boxplus \phi)=L(1,\, \phi)$ since $L(s,\, 1\boxplus \phi)=\zeta(s)L(s,\, \phi)$. Hence we obtain
\begin{align}\label{Main+I(x)}
\sum_{n\geq 1}\lambda_{1\boxplus \phi}(n)W\Big(\frac{n}{X}\Big)=L(1,\, \phi)\tilde{W}(1)X
+\frac{1}{2\pi i}\int_{(-\varepsilon)}\tilde{W}(s)L(s,\, 1\boxplus \phi)X^{s}\d s.
\end{align}
We denote by $I(X)$ the second term with integration on the right hand side of \eqref{Main+I(x)}. Inserting a dyadic smooth partition of unit to the $t$-integral, we get
\begin{align}\label{dyadic}
I(X)=\sum_{T\, \text{dyadic}}I(X, T)
\end{align}
where
\begin{align*}
I(X, T):=\frac{X^{-\varepsilon}}{2\pi}\int_{\mathbb{R}}X^{it}\tilde{W}(-\varepsilon+it)
L(-\varepsilon+it,\, 1\boxplus \phi)V\Big(\frac{t}{T}\Big)\d t
\end{align*}
for some fixed $V$ with compact support.
For $\tilde{W}(s)$, by integration by parts, we have the estimate for any $k\geq 1$
\begin{align}\label{Tilde-W}
\tilde{W}(s)=\frac{(-1)^{k}}{s(s+1)\cdots (s+k-1)}\int_{0}^{\infty}W^{(k)}(u)u^{s+k-1} \d u \ll \frac{1}{|s|^k}\Big(\frac{X}{Y}\Big)^{k-1},
\end{align}
since supp$\, W^{(k)}\in [1/2-Y/X,\, 1/2]\cup [1,\, 1+Y/X]$. This estimate allows us to truncate the $t$-integral of $I(X, T)$ at $T\ll X^{1+\varepsilon}/Y$. In addition, by the upper bounds $L(-\varepsilon+it,\, 1\boxplus \phi)\ll (1+|t|)^{3/2+\varepsilon}$ and \eqref{Tilde-W} with $k=1$, we deduce that
$$
I(X, Y)\ll X^{\varepsilon}T^{3/2+\varepsilon}\ll Y
$$
if $T\ll Y^{2/3-\varepsilon}$. Therefore, by the above arguments, we may impose a constraint $Y^{2/3-\varepsilon}\ll T\ll X^{1+\varepsilon}/Y$ in \eqref{dyadic} with an admissible error term.
We only consider positive $T$'s, since negative $T$'s can be handled similarly.  Next, for $I(X, T)$, by the first equality in \eqref{Tilde-W} with $k=1$, we get
\begin{equation}\label{I(X,T)}
    \begin{aligned}
I(X, T) & =-\frac{X^{-\varepsilon}}{2\pi}\int_{1/3}^{3}W'(u)u^{-\varepsilon}\int_{\mathbb{R}}\frac{(Xu)^{it}}{-\varepsilon+it}
L(-\varepsilon+it,\, 1\boxplus \phi)V\Big(\frac{t}{T}\Big)\d t\, \d u
\\\noalign{\vskip 2,8mm}
& \ll \frac{X^{-\varepsilon}}{T} \sup_{u\in[1/3,\, 3]} \bigg|\int_{\mathbb{R}}(Xu)^{it}
L(-\varepsilon+it,\, 1\boxplus \phi)V\Big(\frac{t}{T}\Big)\d t \bigg|.
\end{aligned}
\end{equation}
Hence, in the following, we only need to consider $J(X, T)$ which is defined by
\begin{align}\label{J(XT)}
J(X, T):=\int_{\mathbb{R}}X^{it}L(-\varepsilon+it,\, 1\boxplus \phi)V\Big(\frac{t}{T}\Big)\d t.
\end{align}
To estimate $J(X, T)$, we shall apply functional equation for $L(-\varepsilon+it,\, 1\boxplus \phi)$ to change the variable $s=-\varepsilon+it$ into $1-s=1+\varepsilon-it$.

By inserting the functional equation \eqref{functional eq} into \eqref{J(XT)}, it follows that
\begin{equation*}
    \begin{aligned}
J(X, T) 
=& \int_{\mathbb{R}}X^{it}\, w\,\bar{\omega}\cdot \Big(\frac{t}{2\pi}\Big)^{3(\frac{1}{2}+\varepsilon)} \Big(\frac{t}{2\pi \text{e}}\Big)^{-3it} L(1+\varepsilon-it,\, 1\boxplus \phi)V\Big(\frac{t}{T}\Big)\d t
\\\noalign{\vskip 2,1mm}
&\ +O\Big(\frac{1}{T}\cdot T^{3/2+\varepsilon}\cdot T\Big)
\\\noalign{\vskip 2,8mm}
\ll&\ T^{3/2+\varepsilon}\bigg|\int_{\mathbb{R}}\ \sum_{n\geq 1}\frac{\lambda_{1\boxplus \phi}(n)}{n^{1+\varepsilon-it}}\,
X^{it} \Big(\frac{t}{2\pi \text{e}}\Big)^{-3it}\, V_1\Big(\frac{t}{T}\Big)\d t\bigg| +T^{3/2+\varepsilon},
\end{aligned}
\end{equation*}
for some smooth compactly supported function $V_1$.

Changing the order of the integral of summation above, and making a change of variable $t=T\xi$, we get
\begin{equation*}
    \begin{aligned}
J(X, T)
\ll&\ T^{5/2+\varepsilon}\bigg|\sum_{n\geq 1}\frac{\lambda_{1\boxplus \phi}(n)}{n^{1+\varepsilon}}
\int_{\mathbb{R}}\ V_1(\xi)\e^{iT\xi\log(nX(\frac{T\xi}{2\pi \e})^{-3})}\d \xi\bigg| +T^{3/2+\varepsilon}.
\end{aligned}
\end{equation*}
For the above integral over $\xi$, by the stationary phase method with
$$
h(\xi)=T\xi\log\Big(\frac{nX}{(\frac{T\xi}{2\pi \e})^3}\Big),\quad \xi_0=\frac{2\pi(nX)^{1/3}}{T},\quad h(\xi_0)=3T\xi_0,\quad h''(\xi_0)=-\frac{3T}{\xi_0}\asymp T,
$$
we get (see e.g. \cite[Proposition 8.2]{BlomerKhanYoung})
\begin{align*}
\int_{\mathbb{R}}\ V_1(\xi)\e^{iT\xi\log(nX(\frac{T\xi}{2\pi \e})^{-3})}\d \xi
=& \frac{\e^{ih(\xi_0)}}{T^{1/2}}W_1(\xi_0)+O\Big(\frac{1}{T^{2020}}\Big)
\\\noalign{\vskip 2,8mm}
=& \frac{\e\big(3 (nX)^{1/3}\big)}{T^{1/2}}W_2\Big(\frac{n}{T^3/X}\Big)+O\Big(\frac{1}{T^{2020}}\Big),
\end{align*}
for some inert functions $W_1, W_2$,
and hence
\begin{equation}\label{J(X,T)-upper}
    \begin{aligned}
J(X, T)\ll& \ T^{2+\varepsilon}\bigg|\sum_{n\geq 1}\frac{\lambda_{1\boxplus \phi}(n)}{n^{1+\varepsilon}}
\, \e\big(3 (nX)^{1/3}\big)W_2\Big(\frac{n}{T^3/X} \Big)\bigg| +T^{3/2+\varepsilon}
\\\noalign{\vskip 2,8mm}
\ll& \ \frac{X^{1+\varepsilon}}{T}\bigg|\sum_{n\geq 1}\lambda_{1\boxplus \phi}(n)
\, \e\big(3 (nX)^{1/3}\big)W_3\Big(\frac{n}{T^3/X} \Big)\bigg| +T^{3/2+\varepsilon},
\end{aligned}
\end{equation}
for some inert function $W_3$.

We may restrict $T$ such that
\begin{align}\label{T}
Y^{2/3-\varepsilon}\ll X^{1/3-\varepsilon}\ll T \ll \frac{X^{1+\varepsilon}}{Y}
\end{align}
since otherwise the above sum over $n$ is empty. Combining \eqref{dyadic}, \eqref{I(X,T)}, \eqref{J(XT)} and \eqref{J(X,T)-upper}, we arrive at
\begin{align}\label{I-upper}
I(X)\ll \, \mathop{\sum_{T\, \text{dyadic}}}_{X^{1/3-\varepsilon}\ll T \ll \frac{X^{1+\varepsilon}}{Y}}\left(\frac{X^{1+\varepsilon}}{T^2}\bigg|\sum_{n\geq 1}\lambda_{1\boxplus \phi}(n)
\, \e\big(3 (nX)^{1/3}\big)W\Big(\frac{n}{T^3/X} \Big)\bigg| +T^{1/2+\varepsilon}\right).
\end{align}
Here $X$ on the right hand side of \eqref{I-upper} should be understood as the original $Xu$ in \eqref{I(X,T)} with  $u\in [1/3,\, 3]$. And $W$ is some smooth function such that supp$\, W\in[1/4,\, 4]$. So we only need to consider the case $n \asymp T^3/X$.

Now we make use of the fact that $\lambda_{1\boxplus \phi}(n)=\sum_{\ell m=n}\lambda_{\phi}(m)$. We can insert a dyadic partition of unit for $\ell$-sum and a dyadic smooth partition of unit for $m$-sum. Hence it suffices to estimate the following sum
$$
B(L, M)=\mathop{\sum_{\ell \sim L}\, \sum_{m\asymp M}}_{LM\asymp \frac{T^3}{X}}\lambda_{\phi}(m)
\, \e\big(3 (\ell mX)^{1/3}\big) V\Big(\frac{m}{M} \Big).
$$
Now we consider two cases. \\
\textit{Case 1}.  $L\gg T^{359/228}X^{-97/152}$. We then rewrite $B(L, M)$ as
\begin{align*}
B(L, M)=\sum_{m\asymp M}\lambda_{\phi}(m)V\Big(\frac{m}{M} \Big)\,
\bigg( \sum_{\substack{\ell \sim L\\ LM\asymp \frac{T^3}{X}}}
\, \e\big(3 (\ell mX)^{1/3}\big) \bigg).
\end{align*}
For the inner sum over $\ell$ we apply the method of exponent pairs with A-process (see for example
\cite[Chapter 3]{GK}), by taking the exponent pair $(p,\, q)$ as
\begin{align*}
\big(p,\, q\big)=\Big(\frac{k}{2k+2},\, \frac{k+h+1}{2k+2}\Big)
=\Big(\frac{13}{194}+\varepsilon,\, \frac{76}{97}+\varepsilon\Big)
\end{align*}
where $(k,\, h)=(\frac{13}{84}+\varepsilon,\, \frac{55}{84}+\varepsilon)$ is also an exponent pair due to Bourgain \cite[Theorem 6]{Bourgain2017decoupling}. Hence we obtain
\begin{equation}\label{case1}
    \begin{aligned}
B(L, M)\ll\, & \sum_{m\asymp M}|\lambda_{\phi}(m)|\,
\Big| \sum_{\substack{\ell \sim L\\ LM\asymp \frac{T^3}{X}}}
\, \e\big(3 (\ell mX)^{1/3}\big) \Big|
\\\noalign{\vskip 1,7mm}
\ll\, & T^{\varepsilon}M\cdot (T/L)^p L^q
\\\noalign{\vskip 1,7mm}
\ll\, & T^{1195/456+\varepsilon}X^{-249/304}.
\end{aligned}
\end{equation}
\textit{Case 2}.  $L\ll T^{359/228}X^{-97/152}$. We rewrite $B(L, M)$ as
\begin{align*}
B(L, M)=\sum_{\ell \sim L}\,
\bigg( \sum_{\substack{m\asymp M\\ LM\asymp \frac{T^3}{X}}}
\, \lambda_{\phi}(m)\e\big(3 (\ell mX)^{1/3}\big)V\Big(\frac{m}{M}\Big)  \bigg).
\end{align*}
For the inner sum over $m$, we employ a result of Jutila (see \cite[Theorem 4.6]{Jutila}).
It is easy to verify that $M^{3/4}\ll T \ll M^{3/2}$
for $L\ll T^{359/228}X^{-97/152}$ and
\begin{align}\label{T2}
X^{165/346} \leq T\leq X^{3/5}.
\end{align}
Therefore
\begin{equation}\label{case2}
    \begin{aligned}
B(L, M)\ll\, & \sum_{\ell \sim L}\,
\Big| \sum_{\substack{m\asymp M\\ LM\asymp \frac{T^3}{X}}}
\, \lambda_{\phi}(m)\e\big(3 (\ell mX)^{1/3}\big)V\Big(\frac{m}{M}\Big) \Big|
\\\noalign{\vskip 1,7mm}
\ll\, & L \cdot T^{1/3}M^{1/2} 
\\\noalign{\vskip 1,7mm}
\ll\, &  T^{1195/456+\varepsilon}X^{-249/304}.
\end{aligned}
\end{equation}

Combining \eqref{I-upper}, \eqref{case1} and \eqref{case2}, we have the following upper bound when $T$ satisfies $X^{165/346} \leq T\leq X^{1/2+\delta+\varepsilon}$:
\begin{equation}\label{I}
    \begin{aligned}
I(X)\ll \, & \frac{X^{1+\varepsilon}}{T^2}\cdot T^{1195/456+\varepsilon}X^{-249/304} + T^{1/2+\varepsilon}
\ll\,  X^{1/2-\delta+\varepsilon},
\end{aligned}
\end{equation}
with $\delta\leq 4/739 \approx 0.00541...$.
If $X^{1/3-\varepsilon} \leq T \leq X^{165/346} $, we use the trivial bound to get
\begin{equation}\label{I2}
    \begin{aligned}
I(X)\ll \, & \frac{X^{1+\varepsilon}}{T^2}\cdot \frac{T^3}{X}  + T^{1/2+\varepsilon}  \ll T X^{\varepsilon} \ll X^{165/346+\varepsilon} \ll X^{1/2-\delta+\varepsilon}.
\end{aligned}
\end{equation}

Finally, putting together the above estimates \eqref{Main+I(x)}, \eqref{T}, \eqref{I} and \eqref{I2}, one can easily complete the proof of Theorem \ref{thm1}.

\begin{remark}
  The exponent pair $(\frac{13}{194}+\varepsilon,\, \frac{76}{97}+\varepsilon)$ is the best known exponent pair we find for our problem. We essentially need to choose an exponent pair $(p,q)$ to minimize $\frac{5+5p-2q}{11+8p-5q}$.
\end{remark}

\vskip 7mm

\bibliographystyle{plain}

\end{document}